\documentclass[12pt]{article}

\usepackage{array}
\usepackage{supertabular}
\usepackage{amsmath,amssymb}
\usepackage{bm}

\newcommand{\xx}{\bm{x}}
\newcommand{\s}{\bm{\varSigma}}
\newcommand{\sss}{\bm{S}}
\newcommand{\cest}{\hat{\bm{\tau}}}
\newcommand{\ceste}[1]{\hat{\tau}_{#1}}
\newcommand{\waj}{\sum_{j=1}^p}
\providecommand{\norm}[1]{\lVert#1\rVert}
\newtheorem{Thm}{Theorem}

\title{Modified estimator of the contribution rates of population eigenvalues}
\author{Yo Sheena\thanks{Department of Economics, Shinshu University}}

\date{April 2010}

\begin{document}
\maketitle

\begin{abstract}
Modified estimators for the contribution rates of population eigenvalues are given under an elliptically contoured distribution. These estimators decrease the bias of the classical estimator, i.e. the sample contribution rates. The improvement of the modified estimators over the classical estimator is proved theoretically in view of their risks. We also checked numerically that the drawback of the classical estimator, namely the underestimation of the dimension in principal component analysis or factor analysis, are corrected in the modification. 
\end{abstract}
\noindent
MSC(2010) {\it  Subject Classification}: Primary 62H12; Secondary 62H25

\section{Introduction}
\label{intro}
Let $\s$ be the population covariance matrix of a $p$-variate random vector $\xx=(x_1,\ldots,x_p)$. Let $\bm{\lambda}=(\lambda_1,\ldots,\lambda_p),\ \lambda_1\geq \cdots \geq \lambda_p$ denote the eigenvalues of $\s$, then the population contribution rates are defined as 
\begin{equation}
\label{def_pop_cont_rate}
\bm{\tau}=(\tau_1,\ldots,\tau_p),\qquad \tau_i=\lambda_i/\sum_{j=1}^p\lambda_j,\quad i=1,\ldots,p .
\end{equation}
The contribution rates play an important role in statistical linear models. Especially in principal component analysis or factor analysis, it gives an important information for the determination of the model's dimension; ``How many principal components substantially represents the total variance ?'' is a basic quercitin in principal component analysis, and the number of factors to be incorporated in a model is a crucial problem in factor analysis.  For this issue, the most simple and widely used methods are the following ones based on the population contribution rates. 
\begin{enumerate}
\item The cumulative percentage of the eigenvalues 
\\
With a cut-off $t^*$, we determine the smallest integer $m$ for which 
\begin{equation}
\label{criterion1}
\sum_{i=1}^m \tau_i \geq t^*
\end{equation}
to be the number of principal components or factors to be retained. Practically a number between $0.7$ and $0.9$ is often chosen as a cut-off $t^*$.
\item The relative size of each eigenvalue
\\
If the $i$th eigenvalue is larger than the average of the population variance $\sum_{\i=1}^p \lambda_i/p$, the corresponding principal component or factor is to be retained. This criteria is equivalent to check whether $\tau_i$ satisfies the inequality

\begin{equation}
\label{criterion2}
\tau_i > p^{-1}. 
\end{equation}

This is also equivalent to ``Kaiser criterion'' in factor analysis, which asserts that the number of the eigenvalues larger than unit of the population {\it correlation} matrix should be the number of factors.
\end{enumerate}

Many methods have been proposed for the determination of dimension relating to principal component analysis or factor analysis (or more generally covariance structure model). See Jolliffe (2002) and Fabrigar et al.(1999), both of which give an extensive review of the methods for choosing a dimension respectively for principal component analysis and exploratory factor analysis. According to Jolliffe (2002)'s classification, there are several categories for the methods other than that based on the population contribution rates; 1) Hypothesis testing method, 2) Information theoretic method , 3) computer-intensive method.

Furthermore it might be better to add another category,  ``large dimensional random matrix method'', if we could name it. For the past decade, while the results have (re)accumulated on ``general asymptotics'' that considers the limiting operation of both $p$ (the dimension) and $n$ (the sample numbers), we have seen much improvement on this method. 
The limiting distribution of the sample eigenvalues under general asymptotics gives some novel ideas for the dimension determination. It is appealing that several simulations show that the arguments based on general asymptotics is effective even if $n$ and $p$ are relatively small.  See e.g. Kritchman and Nadler (2004),   
Ulfarsson and Solo (2008). They propose highly efficient methods for a so-called ``spiked covariance model'' (see the equation \eqref{factor_model}).  For large dimensional random matrix theories used in these papers,  see the references therein. We also refer to Paul (2007), Nadler (2008), Karoui (2009) for more recent developments.  

We should notice that the concept of  ``dimension'' could be rather ambiguous term. In the fields such as physics or chemistry, it is often the case that there exist  ``signals (components)'' and  ``noise'' in its own mechanism with clear distinction. Naturally the analysis of the covariance structure is aimed at ``detecting'' the numbers of the signals (components), as the term ``signal detection'' indicates. On the contrary in psychology or economics, a ``factor''  in its theory is rather abstract object and sometimes impossible to draw the line between the ``factors'' and ``noise''.  There we could only say some factors are trivial while the others are nontrivial. Hence the number of the factors (i.e. dimension) are not considered to preexist before the statistical inference but rather be determined  through the inference so that we can carry out dimension reduction without serious loss of information. We could say it is ``deciding'' the dimension. Considering the both cases, ``detecting'' and ``deciding'' the dimension, it seems that there is no single excellent method that is universally effective. After all we had better choose effective methods according to the purpose of the determination of the dimension and/or the presupposed mechanism of data generation. 

Back to the methods \eqref{criterion1}, \eqref{criterion2} of our concern, their cut-off values seem somewhat ad-hoc without rigorous theoretical background. We naturally raise a question such as  ``Why is $0.9$ for $t^*$ chosen ?'' We only could insist that it is nothing more than conventional criteria for the inference like a given significant level (e.g. 5\%) in a hypothesis test. Nevertheless, they have been widely used and incorporated into many softwares for statistical analysis because of their simplicity and easiness for calculation. They do not suppose any rigid data generation mechanism, which is often suitable for the purpose of  dimension ``decision (reduction)''. We think that the improved inference on the population contribution rates could make some contribution to the better dimension reduction. In this paper we focus ourselves to the point estimation of $\bm{\tau}$ using the sample covariance matrix. 

 Let $\bm{A}$ denote the (unbiased) sample covariance matrix and 
$l^*_1 \geq \cdots \geq l^*_p >0$ be its ordered eigenvalues. Then the sample contribution rates are defined as 
\begin{equation}
\label{def_sample_cont_rate}
d_i=l^*_i/\sum_{j=1}^p l^*_j,\quad 1\leq i \leq p.
\end{equation}
Traditionally (and perhaps almost always) the set of  sample contribution rates has been used for the estimation of the population contribution rates. Hereafter the sample contribution rates as an estimator of $\bm{\tau}$ will be called ``classical estimator'' and denoted by $\cest^{(0)}=(\ceste{1}^{(0)},\ldots,\ceste{p}^{(0)}),$ where $\ceste{i}^{(0)}=d_i,\ i=1,\ldots,p.$

As far as we know, for the estimation of $\bm{\tau}$ there is no other option than the classical estimator. However $\cest^{(0)}$ seems to have nonnegligible bias. It is well known that the sample eigenvalues $\bm{l^*}=(l^*_1,\cdots,l^*_p)$ are biased. Specifically saying, $\lambda_i,\ i=1,\ldots,p$ is majorized by $E(l^*_i),\ i=1,\ldots,p$, that is,
$$
\sum_{j=1}^m \lambda_j \leq \sum_{j=1}^m E(l^*_i),\quad 1\leq m \leq p.
$$
This fact makes us conjecture that $\cest^{(0)}$ is also biased. Let $\sss=(N-1)\bm{A}$, where $N$ is the number of the samples. In the case $\sss$ is distributed as Wishart matrix and $\lambda$'s have no multiplicity, the expected value of $d_i$ is expanded with respect to the degree of freedom $n(=N-1)$ as follows (see the proof in Appendix);
\begin{equation}
\label{assym_expan_E(d_i)}
\begin{split}
E(d_i)&=\tau_i +n^{-1}\left\{2\frac{\lambda_i(\waj \lambda_j^2)}{(\waj \lambda_j)^3}-2\frac{\lambda_i^2}{(\waj \lambda_j)^2}+\frac{\lambda_i}{\waj \lambda_j}\sum_{j\ne i}\frac{\lambda_j}{\lambda_i-\lambda_j}\right\}\\
       &\quad +O(n^{-2})
\end{split}
\end{equation}
The coefficient of the $n^{-1}$ term is complicated, but we easily notice that when $\lambda_i$'s are close to each other, large positive (negative) bias might take place for the smaller (larger) $i$'s.  Note that the similar expansion with respect to both $p$ and $n$ might be possible in view of  ``general asymptotics''. Please refer to Nadler (2008) for a matrix perturbation approach. 

The simulated results of the distribution of $d_i$'s under the condition $\s$ is a identity matrix can be found in Mandel (1972) and Krzanowski (1979). They observe the large bias of  $d_i$'s.  (See also Sugiyama and Tong (1976), Konishi (1977), and Huang and Tseng (1992) for the distribution of $d_i$'s.) Table 1 shows the simulated values of $E(d_i),\ i=1,\ldots,10$ calculated from 10000 random 10-dimensional Wishart matrices with the degree of freedom 30 generated under several patterns of $\bm{\lambda}$. (The total sum of $\lambda_i$'s always equals one, hence $\lambda_i=\tau_i$, $1\leq i \leq p.$ ) As the figures in the table show, it is not seldom that that the bias of the first or last few $E(d_i)$'s surpasses 50\% (sometimes 100\%) of $\lambda_i$, while the sign of  the bias for the middle part of $E(d_i)$'s is quite unstable.
\begin{table}
\caption{Bias of $d_i$} 
\begin{tabular}{|c|c|c|c|c|c|c|c|c|c|c|c|c|c|c|c|c|c|c|c|}
\hline
${\bm{\lambda}_1}$ & ${\bm{\lambda}_2}$ &${\bm{\lambda}_3}$ &${\bm{\lambda}_4}$ &${\bm{\lambda}_5}$ &${\bm{\lambda}_6}$ &${\bm{\lambda}_7}$ &${\bm{\lambda}_8}$ &${\bm{\lambda}_9}$ &${\bm{\lambda}_{10}}$ \\ \hline
$E(d_1)$ & $E(d_2)$ &$E(d_3)$ &$E(d_4)$ &$E(d_5)$ &$E(d_6)$ &$E(d_7)$ &$E(d_8)$ &$E(d_9)$ &$E(d_{10})$ \\ \hline\hline
{\bf 0.100} &{\bf 0.100} &{\bf 0.100} &{\bf 0.100} &{\bf 0.100} &{\bf 0.100} &{\bf 0.100} &{\bf 0.100} &{\bf 0.100} &{\bf 0.100}  \\ \hline
0.216 & 0.172 & 0.142 & 0.118 & 0.097 & 0.080 & 0.064 & 0.050 & 0.037 & 0.024 \\ \hline\hline
{\bf 0.120} &{\bf 0.120} &{\bf 0.120} &{\bf 0.120} &{\bf 0.120} &{\bf 0.080} &{\bf 0.080} &{\bf 0.080} &{\bf 0.080} &{\bf 0.080} \\ \hline
0.225 & 0.176 & 0.143 & 0.117 & 0.095 & 0.077 & 0.061 & 0.047 & 0.035 & 0.023 \\ \hline\hline
{\bf 0.140} &{\bf 0.140} & {\bf 0.140} & {\bf 0.140} & {\bf 0.140} & {\bf 0.060} & {\bf 0.060} & {\bf 0.060} & {\bf 0.060} & {\bf 0.060} \\ \hline
0.245 & 0.189 & 0.149 & 0.116 & 0.090 & 0.069 & 0.053 & 0.040 & 0.029 & 0.019 \\ \hline\hline
{\bf 0.160} &{\bf 0.160} &{\bf 0.160} &{\bf 0.160} &{\bf 0.160} &{\bf 0.040} &{\bf 0.040} &{\bf 0.040} &{\bf 0.040} &{\bf 0.040}  \\ \hline
0.272 & 0.206 & 0.159 & 0.121 & 0.087 & 0.053 & 0.039 & 0.029 & 0.021 & 0.014 \\ \hline\hline
{\bf 0.180}&{\bf 0.180}& {\bf 0.180}& {\bf 0.180}& {\bf 0.180}& {\bf 0.020}& {\bf 0.020}& {\bf 0.020}& {\bf 0.020}& {\bf 0.020} \\ \hline
0.300 & 0.226 & 0.173 & 0.129 & 0.090 & 0.028 & 0.020 & 0.015 & 0.011 & 0.007 \\ \hline\hline
{\bf 0.198} &{\bf 0.198} &{\bf 0.198} &{\bf 0.198} &{\bf 0.198} & {\bf 0.002} &{\bf 0.002} &{\bf 0.002} &{\bf 0.002} &{\bf 0.002} \\ \hline
0.328 & 0.245 & 0.187 & 0.139 & 0.094 & 0.003 & 0.002 & 0.002 & 0.001 & 0.001 \\ \hline\hline
{\bf 0.200} &{\bf 0.089} &{\bf 0.089} &{\bf 0.089} &{\bf 0.089} &{\bf 0.089} &{\bf 0.089} &{\bf 0.089} &{\bf 0.089} &{\bf 0.089}\\ \hline
0.250 & 0.171 & 0.137 & 0.112 & 0.092 & 0.075 & 0.060 & 0.046 & 0.034 & 0.023 \\ \hline\hline
{\bf 0.400} & {\bf 0.067} & {\bf 0.067} & {\bf 0.067} & {\bf 0.067} & {\bf 0.067} & {\bf 0.067} & {\bf 0.067} & {\bf 0.067} & {\bf 0.067}\\ \hline
0.419 & 0.134 & 0.106 & 0.087 & 0.071 & 0.058 & 0.046 & 0.036 & 0.026 & 0.018 \\ \hline\hline
{\bf 0.600} &{\bf 0.044} &{\bf 0.044} & {\bf 0.044} & {\bf 0.044} & {\bf 0.044} & {\bf 0.044} & {\bf 0.044} & {\bf 0.044} & {\bf 0.044} \\ \hline
0.605 & 0.091 & 0.072 & 0.059 & 0.048 & 0.039 & 0.031 & 0.024 & 0.018 & 0.012 \\ \hline\hline
{\bf 0.800} &{\bf 0.022} &{\bf 0.022} &{\bf 0.022} &{\bf 0.022} &{\bf 0.022} &{\bf 0.022} &{\bf 0.022} &{\bf 0.022} &{\bf 0.022}  \\ \hline
0.799 & 0.046 & 0.037 & 0.030 & 0.024 & 0.020 & 0.016 & 0.012 & 0.009 & 0.006 \\ \hline\hline
{\bf 0.990} &{\bf 0.001} &{\bf 0.001} &{\bf 0.001} &{\bf 0.001} &{\bf 0.001} &{\bf 0.001} &{\bf 0.001} &{\bf 0.001} &{\bf 0.001}  \\ \hline
0.990 & 0.002 & 0.002 & 0.002 & 0.001 & 0.001 & 0.001 & 0.001 & 0.000 & 0.000 \\ \hline
\end{tabular}
\end{table}

The aim of this paper is to derive an alternative estimator that modifies the bias of the classical estimator. In the next section, first we show the distribution of the sample contribution rates is identical under a class of elliptical distributions. Second we propose a class of new estimators and show their superiority to the classical estimator under the class of elliptical distributions from a decision theoretic point of view. In the third section, by simulation studies, we clarify other preferable aspects of the new estimator.
%
%
%
%
%
%
%
%
\section{Main Result}
\subsection{Framework}
Before deriving new estimators,  we formulate the estimation problem  of our concern. Let $\bm{x}^{(i)}, i=1,\ldots,N$ be independently and identically distributed  $p$-dimensional sample vectors with a covariance matrix $\s$. Suppose $N>p.$ The sample covariance matrix
$$
\bm{A}=\frac{1}{N-1}\sum_{i=1}^N(\bm{x}^{(i)}-\bar{\bm{x}})(\bm{x}^{(i)}-\bar{\bm{x}})',
$$
where $\bar{\bm{x}}$ is a sample mean vector, is an unbiased estimator of $\s$. We consider the estimation problem of $\bm{\tau}$ (defined by \eqref{def_pop_cont_rate}) based on  $\bm{A}.$ 

We define the following notations;
\begin{align*}
&\bm{X}=\left(
\begin{array}{c}
(\bm{x}^{(1)})' \\
\vdots \\
(\bm{x}^{(N)})'
\end{array}
\right),
\qquad
\bm{Y}=\left(
\begin{array}{c}
(\bm{y}^{(1)})' \\
\vdots \\
(\bm{y}^{(N)})'
\end{array}
\right),\quad \bm{y}^{(i)}=\bm{x}^{(i)}-\bar{\bm{x}},\ i=1,\ldots,N,
\\
&\bm{C}=\bm{I}_N-N^{-1}\bm{1}\bm{1}',
\end{align*}
where $\bm{I}_N$ is the $N$-dimensional identity matrix, and $\bm{1}$ is the $N$-dimensional vector with unit as each element.
We find that
$$
\bm{Y}=\bm{C}\bm{X}.
$$
The expression of $\bm{A}$
$$
\bm{A}=\frac{1}{N-1}\bm{Y}'\bm{Y}
$$
is inconvenient, since the rows of $\bm{Y}$ are linearly constrained. Notice that
$$
\bm{A}=\frac{1}{N-1}\bm{Y}'\bm{Y}=\frac{1}{N-1}\bm{X}'\bm{C}^2\bm{X}=\frac{1}{N-1}\bm{X}'\bm{C}\bm{X}.
$$
Using the decomposition of $\bm{C}=\bm{O_1}\bm{O}_1',\ \bm{O}_1\in V_{N-1,N}$, where $V_{N-1,N}$ is the Stiefel Manifold, if we put $\bm{Z}=\bm{O}'_1\bm{X}$, then we have 
\begin{equation}
\label{gene_form_samp_covar}
\bm{A}=\frac{1}{N-1}\bm{Z}'\bm{Z}.
\end{equation}
The distribution of $\bm{A}$ is determined by $\bm{Z}$ through \eqref{gene_form_samp_covar}, where $\bm{Z}$ is not degenerated.

The most frequently postulated situation is that $\bm{x}_ i\ i=1,\ldots,N$ is independently distributed as the $p$-variate normal distribution $N_p(\bm{\mu},\s)$. Then $(N-1)\bm{A}$ is distributed as Wishart distribution $W_p(N-1,\s)$.
This is distiributionally equivalent to postulating that $\bm{Z}$ is distributed as
$$
\bm{Z}\sim N_{n\times p}(\bm{0}, \bm{I}_n\otimes \s),\quad n=N-1.
$$
The density function of $\bm{Z}$ is proportional to 
\begin{equation}
\label{density_normal}
\exp(-(1/2){\rm tr}\bm{Z}'\bm{Z}\s^{-1})|\s|^{-n/2}.
\end{equation}
One of the natural generalizations of \eqref{density_normal} is an elliptically contoured distribution, the density of  which is given by
\begin{equation}
\label{density_elliptical}
f({\rm tr}\bm{Z}'\bm{Z}\s^{-1})|\s|^{-n/2}.
\end{equation}

We formulate our estimation problem as follows;
\\
$\bm{Z}$ is a $n\times p\ (n\geq p \geq 2)$ random matrix, and its density with respect to Lebesgue measure on $R^{np}$ is given by 
\eqref{density_elliptical} with some function $f(\cdot)$ on $R$, where $\s$ is an unknown positive definite $p$-dimensional matrix. 
We just observe 
\begin{equation}
\label{def_S}
\sss=\bm{Z}'\bm{Z}.
\end{equation}
We consider the estimation of the population contribution rates $\bm{\tau}$ given by \eqref{def_pop_cont_rate} based on $\sss$.
%
%
%
%
%
%
%
%
\subsection{Distribution of the Sample Contribution Rates}
From \eqref{density_elliptical} and \eqref{def_S}, the density of $\sss$ is given by
\begin{equation}
\label{density_S}
c_1 f({\rm tr} \sss\s^{-1}) |\sss|^{(n-p-1)/2} |\s|^{-n/2}
\end{equation}
with some constant $c_1$ (the proof can be found in Appendix). When $f(x)=\exp(-x/2)$ this density function is that of a Wishart distribution. The distribution $\sss$ and the parameter $\s$ are equivariant with respect to the transformations $\sss\to \bm{B}\sss\bm{B}',\ \s \to \bm{B}\s\bm{B}'$ for any $p$-dimensional nonsingular matrix $\bm{B}.$

The eigenvalues of $\sss$ are denoted by $l_i, \ 1\leq i \leq p$. We derive the distribution of the sample contribution rates
$$
d_i=\frac{l_i}{\waj l_j},\quad i=1,\ldots,p.
$$
They are on the hyperplane 
$$
\mathcal{D}=\Bigl\{(d_1,\ldots,d_p)\:\Bigl|\:d_1>\cdots >d_p>0,\quad \sum_{i=1}^pd_i=1\Bigr.\Bigr\}.
$$
We use the coordinate system $\bm{d}=(d_1,\ldots,d_{p-1})$ for $\mathcal{D}$. The range of $\bm{d}$ is given by 
$$
\mathcal{D}^*=\Bigl\{\bm{d}\:\Bigl|\:d_1>\cdots >d_{p-1}>0,\quad \sum_{i=1}^{p-1}d_i<1\Bigr.\Bigr\}.
$$
\begin{Thm}
\label{distributional_result}
Let 
\begin{equation}
\label{spect_decomp_S}
\sss=\bm{H}\bm{L}\bm{H}',\quad \bm{L}={\rm diag}(l_1,\ldots,l_p), \quad \bm{H}\in O(p)
\end{equation}
be the spectral decomposition of $\sss$, where $O(p)$ is the set of $p$-dimensional orthogonal matrices.\\
(i) The density function of  $(\bm{d}, \bm{H})$ with respect to the product measure between Lebesgue measure on $R^{p-1}$ and the invariant probability measure $\mu_{p,p}$ on $O(p)$ is given by
\begin{equation}
\label{final_form_density_d&H}
c_2 |\s|^{-n/2}\: F(\bm{d})\Bigl({\rm tr}\s^{-1}\bm{H}\bm{D}\bm{H}'\Bigr)^{-np/2},
\end{equation}
where $c_2$ is a constant,  $\bm{D}={\rm diag}(d_1,\ldots,d_{p-1},1-\sum_{j=1}^{p-1}d_j),$ and $F(\cdot)$ is a positive-valued function on $\mathcal{D}^*$ that is independent of $f(\cdot)$ in \eqref{density_elliptical} or $\s$.\\
(ii) The density function of $\bm{d}$ with respect to Lebesgue measure on $R^{p-1}$ is given by 
\begin{equation}
\label{final_form_density_d}
c_2\prod_{i=1}^p \tau_i^{-n/2}\: F(\bm{d})\int_{O(p)}\Bigl({\rm tr}\bm{T}^{-1}\bm{H}\bm{D}\bm{H}'\Bigr)^{-np/2} \mu_{p,p}(d\bm{H}),
\end{equation}
where $\bm{T}={\rm diag}(\tau_1,\ldots,\tau_p)$.
\end{Thm}
{\it Proof} \ From \eqref{density_S}, the density function of $\bm{l}=(l_1,\ldots,l_p)$ and $\bm{H}$ with respect to the product measure between Lebesgue measure on $R^p$ and $\mu_{p,p}$ on $O(p)$ is given by 
$$
\tilde{c}\: c_1 |\s|^{-n/2} \prod_{i=1}^p l_i^{(n-p-1)/2} \prod_{i<j}(l_i-l_j) f\Bigl({\rm tr}\s^{-1}\bm{H}\bm{L}\bm{H}'\Bigr)
$$
with some positive constant $\tilde{c}.$ (See e.g. (22) on the p105 of Muirhead(1982).) 

Let $t=\sum_{i=1}^p l_i$ and consider the transformation of coordinates
$$
\bm{l} \to (t,\bm{d}).
$$
Since the Jacorbian is given by $J(\bm{l}\to(t,\bm{d}))=t^{p-1}$, the density function of $t,\bm{d}$ and $\bm{H}$ with respect to $d{\bm{t}}\times d{\bm{d}} \times \mu_{p,p}(d\bm{H})$ is given by 
\begin{equation}
\label{density_t_d_H}
\tilde{c}\:c_1|\s|^{-n/2}\: t^{np/2-1} F(\bm{d})  f\Bigl(t\:{\rm tr}\s^{-1}\bm{H}\bm{D}\bm{H}'\Bigr),
\end{equation}
where 
\begin{align*}
F(\bm{d})&=\Bigl(\prod_{1\leq i<j \leq p-1}(d_i-d_j)\Bigr)\Bigl(\prod_{i=1}^{p-1}(d_i-(1-\sum_{j=1}^{p-1}d_j))\Bigr) \\
&\qquad \times\Bigl(\prod_{i=1}^{p-1}d_i^{(n-p-1)/2}\Bigr)\Bigl(1-\sum_{j=1}^{p-1}d_j\Bigr)^{(n-p-1)/2}.
\end{align*}
Integrate \eqref{density_t_d_H} over $\{t\:|\:0<t<\infty\},$ then 
\begin{align*}
&\tilde{c}\:c_1|\s|^{-n/2}\:  F(\bm{d}) \int_0^\infty t^{np/2-1}f(ta)\: dt \quad \Bigl(a= {\rm tr}\s^{-1}\bm{H}\bm{D}\bm{H}'>0\Bigr)\\
&\ =\tilde{c}\:c_1|\s|^{-n/2}\:  F(\bm{d})\: a^{-np/2}\; \int_0^\infty t^{np/2-1} f(t) dt.
\end{align*}
Now we have the density function of $(\bm{d},\bm{H})$ with respect to $d\bm{d}\times \mu_{p,p}(d\bm{H})$ as
$$
\tilde{c}\:c_1|\s|^{-n/2}\Bigl(\int_0^\infty t^{np/2-1} f(t) dt\Bigr)F(\bm{d})  \Bigl({\rm tr}\s^{-1}\bm{H}\bm{D}\bm{H}'\Bigr)^{-np/2}.
$$
Substituting
$$
\tilde{c}\:c_1\int_0^\infty t^{np/2-1} f(t) dt
$$
with $c_2$, we have \eqref{final_form_density_d&H}. 

Integrating \eqref{final_form_density_d&H} over $O(p)$, we have the density of $\bm{d}$ as
\begin{equation}
\label{density_d}
c_2 \prod_{i=1}^p \lambda_i^{-n/2} F(\bm{d})\int_{O(p)}\Bigl({\rm tr}\s^{-1}\bm{H}\bm{D}\bm{H}'\Bigr)^{-np/2}\mu_{p,p}(d\bm{H}).
\end{equation}
Let 
$$
\s=\tilde{\bm{H}}'\bm{\varLambda}\tilde{\bm{H}},\quad \bm{\varLambda}={\rm diag}(\lambda_1,\ldots,\lambda_p),\quad \tilde{\bm{H}}\in  O(p)
$$
be the spectral decomposition  of $\s.$ Since $\mu_{p,p}$ is the invariant probability on $O(p)$, $\tilde{\bm{H}}\bm{H}$ has the same distribution as $\bm{H}$. Therefore \eqref{density_d} equals
$$
c_2 \prod_{i=1}^p \lambda_i^{-n/2} F(\bm{d})\int_{O(p)}\Bigl({\rm tr}\bm{\varLambda}^{-1}\bm{H}\bm{D}\bm{H}'\Bigr)^{-np/2}\mu_{p,p}(d\bm{H}).
$$
Substituting $\lambda_i$ with $\tau_i\waj \lambda_j$, we have \eqref{final_form_density_d}. \hfill \rule{5pt}{10pt}
\\

It is noteworthy that the distribution of $\bm{d}$ is independent of $f(\cdot)$ in \eqref{density_elliptical} and depends on $\s$ only through the population contribution rates $\bm{\tau}$.  Note that Johnson and Grayvill (1972) deals with the distribution of the sample contribution rates when $\s=I_p$.
\subsection{New Estimator}
In order to derive a new estimator $\hat{\bm{\tau}}(\bm{d})=(\hat{\tau}_1(\bm{d}),\ldots,\hat{\tau}_p(\bm{d}))$ of $\bm{\tau}=(\tau_1,\ldots,\tau_p)$ that has a certain superiority to the classical estimator, we take a decision theoretic approach here, that is,  we compare estimators via their risks with respect to a certain loss function. Straightforward approach is to use a loss function that directly measures the distance between $\bm{\tau}$ and its estimator $\hat{\bm{\tau}}$. However the exact property of the sample contribution rates $\bm{d}$ are difficult to derive under the assumption of small samples. Instead we evaluate the performance of $\hat{\bm{\tau}}(\bm{d})$ as the components of an estimator of $\s$. 

Consider an estimator of $\s$ combining $\hat{\bm{\tau}}$ and the sample eigenvectors $\bm{H}$ in \eqref{spect_decomp_S} as follows;
\begin{equation}
\label{combined_estimator_S}
\hat\s=\bm{H} \hat{\bm{T}}(\bm{d})\bm{H}',\quad \hat{\bm{T}}(\bm{d})={\rm diag}(\hat\tau_1(\bm{d}),\ldots,\hat\tau_p(\bm{d})).
\end{equation}
Our approach is based on the following observations. According to Theorem \ref {distributional_result}, the distribution of $\bm{d}$ is determined by $\bm{\tau}$. Therefore we can suppose that ${\rm tr \s}=\sum_{i=1}^p \lambda_i=1$ without loss of generality. In this case, the population contribution rates are equal to the population eigenvalues. In addition, the sample eigenvectors ${\bf H}$ are M.L.E. , hence consistent under the large sample asymptotics if $f(\cdot)$ in (10) is monotonically decreasing. (See Paul (2007) and Nadler (2008) for the discrepancy between the sample eigenvectors and the population counterparts for a large-dimensional matrix.) Therefore \eqref{combined_estimator_S} is supposed to be a good estimator of $\s$ if  $\hat{\bm{\tau}}$ is a good estimator of $\bm{\tau}.$

The most common loss function about $\s$ and $\hat\s$ is the entropy loss function (Stein's loss function) 
\begin{equation}
\label{entropy_loss}
L(\hat\s,\s)={\rm tr}(\hat{\s}\s^{-1})-\log(|\hat{\s}\s^{-1}|)-p.
\end{equation}
We evaluate the performance of $\hat\s$ through its risk with respect to this loss function.
\\

We consider one class of simple estimators given by 
\begin{equation}
\label{est_tau_*}
\hat{\bm{\tau}}^*=(\hat{\tau}_1^*,\ldots,\hat{\tau}_p^*),\quad \hat{\tau}_i^*=\beta_i^* d_i,\ 1\leq i \leq p,
\end{equation}
where $\beta_i^*,\ i=1,\ldots,p$ are positive constants. The classical estimator denoted by $\hat{\bm{\tau}}^0$ is given by
\begin{equation}
\label{est_tau_0}
\hat{\bm{\tau}}^0=(\hat{\tau}_1^0,\ldots,\hat{\tau}_p^0),\quad \hat{\tau}_i^0=\beta_i^0 d_i,\ 1\leq i \leq  p,
\end{equation}
where $\beta_i^0 =1,\ i=1,\ldots,p$. Correspondingly we define the two estimators $\hat{\s}^*$ and $\hat{\s}^0$ as follows;
\begin{eqnarray}
\label{est_sigma_*}
\hat{\s}^*&=\bm{H}\hat{\bm{T}}^*\bm{H}',\quad \hat{\bm{T}}^*={\rm diag}(\hat{\tau}_1^*,\ldots,\hat{\tau}_p^*),\\
\label{est_sigma_0}
\hat{\s}^0&=\bm{H}\hat{\bm{T}}^0\bm{H}',\quad \hat{\bm{T}}^0={\rm diag}(\hat{\tau}_1^0,\ldots,\hat{\tau}_p^0).
\end{eqnarray}

We have the following result on the superiority of $\hat{\s}^*$ to $\hat{\s}^0$.
\begin{Thm}
\label{dominance_result}
If $\beta_i^*\ ( i=1,\ldots,p)$ satisfy the following three conditions, then $\hat{\s}^*$ dominates $\hat{\s}^0$ with respect to the loss function \eqref{entropy_loss}.\\
 For some $m$ $(1\leq m \leq p-1)$, the next two  inequalities hold;
\begin{align}
\label{Condition 1}
&0< \beta_1^* \leq \cdots \leq \beta_m^* \leq 1 \leq \beta_{m+1}^*\leq \cdots \leq \beta_p^*,  \\
\label{Condition 2}
&\sum_{i=1}^m (n+p-1-2i)(\beta_i^*-1)+\sum_{i=m+1}^p(n+p+1-2i)(\beta_i^*-1)\leq 0.
\end{align}
Moreover, the third inequality
\begin{equation}
\label{Condition 3}
\sum_{i=1}^p (\beta_i^*)^{-1} \leq p
\end{equation}
holds.
\end{Thm}
{\it Proof} \ Since both $\hat{\s}^*$ and $\hat{\s}^0$ is the function of  $\bm{d}$ and $\bm{H}$, their distributions are independent of $f(\cdot)$ from the result (i) of Theorem \ref{distributional_result}. Therefore  we can suppose $f(x)=\exp(-x/2)$, that is, $\sss$ is distributed as a Wishart matrix;  
\begin{equation}
\label{s_Wishart}
\sss \sim W_p(n,\s).
\end{equation}

If $\sss$ is distributed as in \eqref{s_Wishart}, the following Stein-Haff identity holds. (Exactly speaking, it is the application of Stein-Haff identity to an orthogonally equivariant estimator, see e.g. Lemma 2.1 of Dey and Srinivasan (1986));
\\
Suppose $\sss$ is decomposed as in \eqref{spect_decomp_S} and $\hat{\s}$ is given by
$$
\hat{\s}=\bm{H}{\rm diag}(\phi_1(\bm{l}),\ldots, \phi_p(\bm{l}))\bm{H}', \quad \bm{l}=(l_1,\ldots,l_p).
$$
Then 
$$
E[{\rm tr}(\hat{\s}\s^{-1})]=E[G(\hat\s,\bm{l})],
$$
where
\begin{equation}
\label{form_G}
G(\hat\s,\bm{l})=2\sum_{1\leq i < j \leq p} \frac{\phi_i(\bm{l})-\phi_j(\bm{l})}{l_i-l_j} +2\sum_{i=1}^p \frac{\partial \phi_i(\bm{l})}{\partial l_i}
+(n-p-1) \sum_{i=1}^p\frac{\phi_i(\bm{l})}{l_i}.
\end{equation}
If we use this identity, we have the following equation.
\begin{equation}
\label{transform_risk_difference}
\begin{split}
&E[L(\hat{\s}^*,\s)]-E[L(\hat{\s}^0,\s)]\\
&=E[{\rm tr}(\hat{\s}^*\s^{-1})-\log|\hat{\s}^*|]-E[{\rm tr}(\hat{\s}^0\s^{-1})-\log|\hat{\s}^0|]\\
&=E[G(\hat{\s}^*,\bm{l})-G(\hat{\s}^0,\bm{l})-\log|\hat{\s}^*|+\log|\hat{\s}^0|].
\end{split}
\end{equation}
Substituting \eqref{est_tau_*} and \eqref{est_sigma_*} into \eqref{form_G}, followed by simple calculation, we have
\begin{equation}
\label{specific_form_G}
\begin{split}
G(\hat{\s}^*,\bm{l})&=\frac{2}{\waj l_j} \sum_{1\leq i<j \leq p}\frac{\beta_i^* l_i-\beta_j^* l_j}{l_i-l_j} +\frac{2}{\waj l_j}\sum_{i=1}^p \beta_i^* \frac{\waj l_j-l_i}{\waj l_j}\\
&\quad+\frac{n-p-1}{\waj l_j} \sum_{i=1}^p \beta_i^* \\
&=\frac{1}{\waj l_j}\Bigl\{2\sum_{1\leq i <j \leq p}\frac{\beta_i^*l_i-\beta_i^*l_j+\beta_i^* l_j-\beta_j^*l_j}{l_i-l_j} \\
&\quad+2\sum_{i=1}^p\beta_i^*(1-d_i)+(n-p-1)\sum_{i=1}^p\beta_i^*\Bigr\}\\
&=\frac{1}{\waj l_j}\Bigl\{2\sum_{1\leq i <j \leq p}(\beta_i^*-\beta_j^*)\frac{l_j}{l_i-l_j}+2\sum_{i=1}^p(p-i)\beta_i^*\\
&\quad+2\sum_{i=1}^p\beta_i^*(1-d_i)+(n-p-1)\sum_{i=1}^p\beta_i^*\Bigr\}\\
&=\frac{1}{\waj l_j}\Bigl\{2\sum_{1\leq i <j \leq p}(\beta_i^*-\beta_j^*)\frac{l_j}{l_i-l_j}+\sum_{i=1}^p(n+p+1-2i-2d_i)\beta_i^*\Bigr\}
\end{split}
\end{equation}
If we substitute $\beta_i^*$ in \eqref{specific_form_G} with $\beta_i^0=1$ $(1\leq i \leq p)$, we have $G(\hat{\s}^0,\bm{l})$. From these results, the inside of the brackets of the right-hand side in \eqref{transform_risk_difference} turns out to be 
\begin{equation}
\label{within_bracket_risk_diff}
\begin{split}
&\Bigl(\waj l_j\Bigr)^{-1} \Bigl\{2\sum_{1\leq i <j \leq p}(\beta_i^*-\beta_j^*)\frac{l_j}{l_i-l_j}+\sum_{i=1}^p(n+p+1-2i-2d_i)(\beta_i^*-1)\Bigr\}\\
&\quad-\log\prod_{i=1}^p\beta_i^*.
\end{split}
\end{equation}
Since $\beta_i^* \leq \beta_j^*$ for $1\leq i<j \leq p$ from \eqref{Condition 1}, \eqref{within_bracket_risk_diff} is less than or equal to
\begin{equation}
\label{within_bracket_risk_diff_2}
\Bigl(\waj l_j\Bigr)^{-1} \sum_{i=1}^p(n+p+1-2i-2d_i)(\beta_i^*-1)-\log\prod_{i=1}^p\beta_i^*. 
\end{equation}
\eqref{Condition 1} says that $\beta_i^*\leq 1$ if $1\leq i \leq m$ and that $\beta_i^*\geq 1$ if $m+1\leq i \leq p$. Using this fact together with the inequality $0\leq d_i \leq 1,\ \forall i$, we notice that \eqref{within_bracket_risk_diff_2} is less than or equal to
$$
\Bigl(\waj l_j\Bigr)^{-1}\Bigl\{ \sum_{i=1}^m(n+p-1-2i)(\beta_i^*-1)+\sum_{i=m+1}^p(n+p+1-2i)(\beta_i^*-1)\Bigr\}-\log\prod_{i=1}^p\beta_i^* .
$$
From \eqref{Condition 2}, this is less than or equal to $\sum_{i=1}^p \log (\beta_i^*)^{-1}$. Because of the inequality $\log (x+1) \leq x,\ \forall x>-1$, we have
$$
\sum_{i=1}^p \log (\beta_i^*)^{-1}=\sum_{i=1}^p \log \{(\beta_i^*)^{-1}-1+1\} \leq \sum_{i=1}^p \{(\beta_i^*)^{-1}-1\},
$$
which is nonpositive by \eqref{Condition 3}. \hfill \rule{5pt}{10pt}
\\

\eqref{Condition 1} of Theorem \ref{dominance_result} means $\bm{\tau}^*$ modifies the bias of the classical estimator which we mentioned in Section \ref{intro}, since lighter weight is given to $d_i$ for the smaller $i$'s and heavier weight for the larger $i$'s.

Choose an  integer $q$ such that $1\leq q \leq p/2-1.$ Let $\beta_i^{(q)}\ (i=1,\ldots,p)$ be defined as 
\begin{equation}
\label{estimator q}
\beta_i^{(q)}=\left\{
\begin{split}
&n(n+p-2q+1-2i)^{-1} &&\text{if $1\leq i \leq m-q$,}\\
&\qquad \quad 1                    &&\text{if $m-q+1 \leq i  \leq p-m+q$,}\\
&n(n+p+2q+1-2i)^{-1}  &&\text{if $p-m+q+1 \leq i \leq p$,}
\end{split}
\right.
\end{equation}
where $m=[p/2]$, i.e., the largest integer that does not exceed $p/2.$ Then $\beta_i^{(q)}\ (i=1,\ldots,p)$ satisfy \eqref{Condition 1}, \eqref{Condition 2} and \eqref{Condition 3}. In fact, \eqref{Condition 1} is clearly  satisfied from the definition. \eqref{Condition 2} and \eqref{Condition 3} are also satisfied as follows;
\begin{align*}
&\sum_{i=1}^m (n+p-1-2i)(\beta_i^{(q)}-1)+\sum_{i=m+1}^p (n+p+1-2i)(\beta_i^{(q)}-1)\\
&=\sum_{i=1}^{m-q} (n+p-1-2i)(\beta_i^{(q)}-1)+\sum_{i=p-m+q+1}^p (n+p+1-2i)(\beta_i^{(q)}-1)\\
&\leq \sum_{i=1}^{m-q} (n+p-2q+1-2i)(\beta_i^{(q)}-1)+\sum_{i=p-m+q+1}^p (n+p+2q+1-2i)(\beta_i^{(q)}-1)\\
&=\sum_{i=1}^{m-q} (-p+2q-1+2i)+\sum_{i=p-m+q+1}^p (-p-2q-1+2i)\\
&=(m-q)\{(-p+2q-1)+(m-q+1)+(-p-2q-1)+(2p-m+q+1)\}\\
&=0,
\end{align*}
\begin{align*}
&\sum_{i=1}^p(\beta_i^{(q)})^{-1}-p\\
&=\sum_{i=1}^p\{(\beta_i^{(q)})^{-1}-1\}\\
&=n^{-1}\Bigl\{\sum_{i=1}^{m-q} (p-2q+1-2i)+\sum_{i=p-m+q+1}^p (p+2q+1-2i)\Bigr\}\\
&=0.
\end{align*}
We give two examples of the estimators that satisfy the three conditions in Theorem \ref{dominance_result}. Let $q=1$, then
\begin{equation}
\label{estimator 1}
\beta_i^{(1)}=\left\{
\begin{split}
&n(n+p-1-2i)^{-1} &&\text{if $1\leq i \leq m-1$,}\\
&\qquad \quad 1                    &&\text{if $m \leq i  \leq p-m+1$,}\\
&n(n+p+3-2i)^{-1}  &&\text{if $p-m+2 \leq i \leq p$,}
\end{split}
\right.
\end{equation}
The specific value of $\beta_i^{(1)}\:(1\leq i \leq p)$ is given as follows; \\
if $p$ is even
\\
\\
\begin{tabular}{c|c|c|c|c|c|c|c|c|c}
$\beta_1^{(1)}$   & $\beta_2^{(1)}$  & $\cdots$ & $\beta_{m-1}^{(1)}$ &$\beta_{m}^{(1)}$ &$\beta_{m+1}^{(1)}$ & $\beta_{m+2}^{(1)}$ & $\cdots$ &
$\beta_{p-1}^{(1)}$ & $\beta_{p}^{(1)}$ \\ \hline
$\frac{n}{n+p-3}$ & $\frac{n}{n+p-5}$ & $\cdots $ &$\frac{n}{n+1}$ & $1$ & $1$ & $\frac{n}{n-1}$ & $\cdots $ & $\frac{n}{n-p+5}$ & $\frac{n}{n-p+3}$ \\
\end{tabular}
\\
\\
\\
if $p$ is odd
\\
\\
\begin{tabular}{c|c|c|c|c|c|c|c|c|c|c}
$\beta_1^{(1)}$   & $\beta_2^{(1)}$  & $\cdots$ & $\beta_{m-1}^{(1)}$ &$\beta_{m}^{(1)}$ &$\beta_{m+1}^{(1)}$ & $\beta_{m+2}^{(1)}$ &$\beta_{m+3}^{(1)}$ & $\cdots$ &
$\beta_{p-1}^{(1)}$ & $\beta_{p}^{(1)}$ \\ \hline
$\frac{n}{n+p-3}$ & $\frac{n}{n+p-5}$ & $\cdots $ &$\frac{n}{n+2}$ & $1$ &$1$&  $1$ & $\frac{n}{n-2}$ & $\cdots $ & $\frac{n}{n-p+5}$ & $\frac{n}{n-p+3}$ \\
\end{tabular}
\\
\\
\\
The estimator 
\begin{equation}
\hat{\s}^{(1)}=\bm{H}\hat{\bm{T}}^{(1)}\bm{H}',\quad \hat{\bm{T}}^{(1)}={\rm diag}(\hat{\tau}_1^{(1)},\ldots,\hat{\tau}_p^{(1)}),
\end{equation}
where $\hat{\tau}_i^{(1)}=\beta_i^{(1)} d_i,\ 1\leq i \leq p$, dominates $\hat{\s}^0$ if $p\geq 4.$ 
Another estimator that satisfies the three conditions of Theorem \ref{dominance_result} is given by $q=2$, which leads to 
\begin{equation}
\label{estimator 2}
\beta_i^{(2)}=\left\{
\begin{split}
&n(n+p-3-2i)^{-1} &&\text{if $1\leq i \leq m-2$,}\\
&\qquad \quad 1                    &&\text{if $m-1 \leq i  \leq p-m+2$,}\\
&n(n+p+5-2i)^{-1}  &&\text{if $p-m+3 \leq i \leq p$.}
\end{split}
\right.
\end{equation}
The specific value of $\beta_i^{(2)}\:(1\leq i \leq p)$ is given as follows; \\
if $p$ is even
\\
\\
\begin{tabular}{c|c|c|c|c|c|c|c|c|c|c}
$\beta_1^{(2)}$   & $\beta_2^{(2)}$  & $\cdots$ & $\beta_{m-2}^{(2)}$ &$\beta_{m-1}^{(2)}$ & $\cdots$ &$\beta_{m+2}^{(2)}$ & $\beta_{m+3}^{(2)}$ & $\cdots$ &
$\beta_{p-1}^{(2)}$ & $\beta_{p}^{(2)}$ \\ \hline
$\frac{n}{n+p-5}$ & $\frac{n}{n+p-7}$ & $\cdots $ &$\frac{n}{n+1}$ & $1$ & $\cdots$ & $1$ & $\frac{n}{n-1}$ & $\cdots $ & $\frac{n}{n-p+7}$ & $\frac{n}{n-p+5}$ \\
\end{tabular}
\\
\\
\\
if $p$ is odd
\\
\\
\begin{tabular}{c|c|c|c|c|c|c|c|c|c|c}
$\beta_1^{(2)}$   & $\beta_2^{(2)}$  & $\cdots$ & $\beta_{m-2}^{(1)}$ &$\beta_{m-1}^{(2)}$& $\cdots$ &$\beta_{m+3}^{(2)}$ & $\beta_{m+4}^{(2)}$ & $\cdots$ &
$\beta_{p-1}^{(2)}$ & $\beta_{p}^{(2)}$ \\ \hline
$\frac{n}{n+p-5}$ & $\frac{n}{n+p-7}$ & $\cdots $ &$\frac{n}{n+2}$ & $1$ & $\cdots$ &$1$& $\frac{n}{n-2}$ & $\cdots $ & $\frac{n}{n-p+7}$ & $\frac{n}{n-p+5}$ \\
\end{tabular}
\\
\\
\\
Note that 
\begin{equation}
\label{relation_tau1_tau2}
\begin{split}
\beta_i^{(1)} \leq \beta_i^{(2)} \leq \beta_i^{(0)}(\equiv 1) &\text{ for $1\leq i \leq m$ },\\
\beta_i^{(1)} \geq \beta_i^{(2)} \geq \beta_i^{(0)}(\equiv 1)&\text{ for $m+1\leq i \leq p.$}
\end{split}
\end{equation}
The estimator 
\begin{equation}
\hat{\s}^{(2)}=\bm{H}\hat{\bm{T}}^{(2)}\bm{H}',\quad \hat{\bm{T}}^{(2)}={\rm diag}(\hat{\tau}_1^{(2)},\ldots,\hat{\tau}_p^{(2)}),
\end{equation}
where $\hat{\tau}_i^{(2)}=\beta_i^{(2)} d_i,\ 1\leq i \leq p$, dominates $\hat{\s}^0$ if $p\geq 6.$ 
\section{Simulation Study}
In this section, we examine by simulation other preferable properties of  the new estimator $\hat{\bm{\tau}}^*=(\hat{\tau}_1^*,\ldots,\hat{\tau}_p^*),$
$$
\hat{\tau}_i^*=\beta_i^* d_i,\ 1\leq i \leq p,
$$
where $\beta_i^*$'s satisfy \eqref{Condition 1}--\eqref{Condition 3}, especially when  $\beta_i^*=\beta_i^{(1)}$ or $\beta_i^{(2)} ,$$\;\ 1\leq i \leq p.$
We use the notation $\hat{\bm{\tau}}^{(j)}=(\hat{\tau}_1^{(j)},\ldots,\hat{\tau}_p^{(j)}),$
$$
\hat{\tau}_i^{(j)}=\beta_i^{(j)} d_i, \qquad j=1,2,\quad 1\leq i \leq p.
$$
\subsection{Risk Comparison}
We will compare the estimator $\hat{\bm{\tau}}^{(1)}$ and $\hat{\bm{\tau}}^{(2)}$ with the classical estimator $\hat{\bm{\tau}}^{(0)}$ through their risks with respect to the quadratic loss function;
\begin{equation}
\label{ql}
QL(\hat{\bm{\tau}},\bm{\tau})=\sum_{i=1}^p (1-\hat{\tau}_i/\tau_i)^2.
\end{equation}
According to Theorem \ref{dominance_result}, the plug-in estimator $\hat{\s}^{(j)}$ made from  $\hat{\bm{\tau}}^{(j)}$ $j=1,2$ dominates another plug-in estimator $\hat{\s}^{(0)}$ from  $\hat{\bm{\tau}}^{(0)}$ with respect to the entropy loss function. We are interested in a more direct comparison among $\hat{\bm{\tau}}^{(0)}$, $\hat{\bm{\tau}}^{(1)}$ and $\hat{\bm{\tau}}^{(2)}$ using \eqref{ql}.

We generated 10000 random 10-dimensional Wishart matrices with the degree of freedom 30 under several patterns of the population contribution rates, $\bm{\tau}=(\tau_1,\ldots, \tau_{10})$. The Table \ref{quadratic_risk} shows the simulation result, where the first 10 numbers in each row are the population contribution rates and the last three numbers are the simulated risks for the three estimators $\cest^{(j)},\:j=0,1,2$ (all the numbers are rounded to the second decimal place).
The risk of $\hat{\bm{\tau}}^{(1)}$ is smaller than that of $\hat{\bm{\tau}}^{(0)}$ by 30\% to 40\%. Since $\hat{\bm{\tau}}^{(2)}$ is located between  $\hat{\bm{\tau}}^{(0)}$ and  $\hat{\bm{\tau}}^{(1)}$ (see \eqref{relation_tau1_tau2}), its risk reduction is smaller than $\hat{\bm{\tau}}^{(1)}$. Nevertheless it still reduces the risk by 17\% to 30\% compared to $\hat{\bm{\tau}}^{(0)}$. From these results, we can conclude that the new estimators are substantially improved over the classical estimator in view of the quadratic risk.
 
\begin{table}
\caption{Risk w.r.t. Quadratic Loss}
\label{quadratic_risk}
\begin{tabular}{|c|c|c|c|c|c|c|c|c|c|c|c|c|}
\hline
$\tau_1$ & $\tau_2$ & $\tau_3$ & $\tau_4$ & $\tau_5$ & $\tau_6$ & $\tau_7$ & $\tau_8$ & $\tau_9$ & $\tau_{10}$
 & $\hat{\bm{\tau}}^{(0)}$ & $\hat{\bm{\tau}}^{(1)}$ & $\hat{\bm{\tau}}^{(2)}$\\ \hline
0.10 & 0.10 & 0.10 & 0.10 & 0.10 & 0.10 & 0.10 & 0.10 & 0.10 & 0.10 & 3.57 & 2.11 & 2.56 \\ \hline
0.11 & 0.11 & 0.11 & 0.11 & 0.11 & 0.09 & 0.09 & 0.09 & 0.09 & 0.09 & 2.71 & 1.57 & 1.91 \\ \hline
0.12 & 0.12 & 0.12 & 0.12 & 0.12 & 0.08 & 0.08 & 0.08 & 0.08 & 0.08 & 2.20 & 1.25 & 1.53 \\ \hline
0.13 & 0.13 & 0.13 & 0.13 & 0.13 & 0.07 & 0.07 & 0.07 & 0.07 & 0.07 & 1.96 & 1.13 & 1.36 \\ \hline
0.14 & 0.14 & 0.14 & 0.14 & 0.14 & 0.06 & 0.06 & 0.06 & 0.06 & 0.06 & 1.85 & 1.10 & 1.31 \\ \hline
0.15 & 0.15 & 0.15 & 0.15 & 0.15 & 0.05 & 0.05 & 0.05 & 0.05 & 0.05 & 1.84 & 1.14 & 1.32 \\ \hline
0.16 & 0.16 & 0.16 & 0.16 & 0.16 & 0.04 & 0.04 & 0.04 & 0.04 & 0.04 & 1.87 & 1.21 & 1.38 \\ \hline
0.17 & 0.17 & 0.17 & 0.17 & 0.17 & 0.03 & 0.03 & 0.03 & 0.03 & 0.03 & 1.89 & 1.27 & 1.42 \\ \hline
0.18 & 0.18 & 0.18 & 0.18 & 0.18 & 0.02 & 0.02 & 0.02 & 0.02 & 0.02 & 1.91 & 1.32 & 1.46 \\ \hline
0.19 & 0.19 & 0.19 & 0.19 & 0.19 & 0.01 & 0.01 & 0.01 & 0.01 & 0.01 & 1.94 & 1.37 & 1.50 \\ \hline
0.20 & 0.09 & 0.09 & 0.09 & 0.09 & 0.09 & 0.09 & 0.09 & 0.09 & 0.09 & 2.71 & 1.78 & 2.15 \\ \hline
0.30 & 0.08 & 0.08 & 0.08 & 0.08 & 0.08 & 0.08 & 0.08 & 0.08 & 0.08 & 2.86 & 1.92 & 2.31 \\ \hline
0.40 & 0.07 & 0.07 & 0.07 & 0.07 & 0.07 & 0.07 & 0.07 & 0.07 & 0.07 & 2.97 & 2.02 & 2.42 \\ \hline
0.50 & 0.06 & 0.06 & 0.06 & 0.06 & 0.06 & 0.06 & 0.06 & 0.06 & 0.06 & 3.09 & 2.12 & 2.53 \\ \hline
0.60 & 0.04 & 0.04 & 0.04 & 0.04 & 0.04 & 0.04 & 0.04 & 0.04 & 0.04 & 3.24 & 2.23 & 2.66 \\ \hline
0.70 & 0.03 & 0.03 & 0.03 & 0.03 & 0.03 & 0.03 & 0.03 & 0.03 & 0.03 & 3.44 & 2.39 & 2.84 \\ \hline
0.80 & 0.02 & 0.02 & 0.02 & 0.02 & 0.02 & 0.02 & 0.02 & 0.02 & 0.02 & 3.65 & 2.56 & 3.03 \\ \hline
0.90 & 0.01 & 0.01 & 0.01 & 0.01 & 0.01 & 0.01 & 0.01 & 0.01 & 0.01 & 3.94 & 2.79 & 3.28 \\ \hline
\end{tabular}
\end{table}
\subsection{Estimation of Dimension}
As we mentioned in Section \ref{intro}, $\bm{\tau}$, the population contribution rates, is one of  the most basic tools for deciding the dimension in principal component analysis or factor analysis. As the first step in deciding the dimension, the choice of an estimator $\cest$ for $\bm{\tau}$ is an important task, hence we are interested in how the new estimator, $\hat{\bm{\tau}}^{(*)}$, makes a difference compared to the classical estimator, $\hat{\bm{\tau}}^{(0)}$, in the decision of the dimension.

Suppose that $\bm{x}=(x_1,\ldots, x_p)'$ are generated in the following $m$-factor model ;
\begin{equation}
\label{factor_model}
\bm{x}=\bm{a}+\bm{B}\bm{z}+\bm{e},
\end{equation}
where $\bm{a}$ is a constant $p$-dimensional vector, $\bm{B}$ is a $p \times m\; (p\geq m)$ factor loading constant matrix with the rank of $m$, $\bm{z}$ is a $m$-dimensional random factor, and $\bm{e}$ is the $p$-dimensional error term which is independent of $\bm{z}$. If we suppose the covariance matrices of $\bm{z}$ and $\bm{e}$ are respectively given by
$$
V(\bm{z})=\bm{\varSigma}_0,\quad V(\bm{e})=\sigma^2 \bm{I}_p,
$$
then the covariance matrix $\s$ of $\bm{x}$ equals 
\begin{equation}
\label{compounded_Sigma}
\s=\bm{B}\bm{\varSigma}_0\bm{B}'+\sigma^2 \bm{I}_p.
\end{equation}
If we denote the eigenvalues of $\bm{B}\bm{\varSigma}_0\bm{B}'$ by $\xi_i\;(1\leq i \leq p)$, then $\lambda_i\; (1\leq i \leq p)$, the eigenvalues of $\s$, are given by
\begin{equation}
\label{compounded_lambda}
\lambda_i=\left\{
\begin{split}
&\xi_i+\sigma^2 &&\text{ if $i=1,\ldots, m,$}\\
&\sigma^2 &&\text{ if $i=m+1,\ldots,p,$}
\end{split}
\right.
\end{equation}
since $\xi_i=0,\:(m+1\leq i \leq p).$ In the large-dimensional random matrix theory, the model \eqref{factor_model} is called ``spiked covariance model'', since the part of $\xi_i+\sigma^2$ seems spiked into the long flat part $\sigma^2$. If $\bm{z}$ and $\bm{e}$ have normal distributions, then $\bm{x}$ is also normally distributed with the covariance matrix \eqref{compounded_Sigma}, hence we can suppose that $\sss\sim W_p(n,\s)$.

We made a simulation under the condition
$$
\sss \sim W_{10}(n,\bm{\Lambda}),\quad n=30, 100,\quad \bm{\Lambda}={\rm diag}(\lambda_1,\ldots,\lambda_p),
$$
where $\lambda_i\;(1\leq i \leq p)$ is given by \eqref{compounded_lambda} with ten patterns of $(m,\xi_1,\ldots,\xi_m)$ and $\sigma^2$ fixed to be unit. $\sss$ is generated 10000 times, and for each time we recorded the dimension decided from the six methods composed by the combination of two criterions in Section\ref{intro} and three estimators $\hat{\bm{\tau}}^{(0)}$, $\hat{\bm{\tau}}^{(1)}$, $\hat{\bm{\tau}}^{(2)}$; the two criterions are ``the cumulative percentage of the (estimated) population eigenvalues'' with $t^*=0.8$ in \eqref{criterion1} (say criterion 1) and ``the relative size of each (estimated) population eigenvalue'' (say criterion 2).

\begin{table}
\caption{Histogram of Estimated Dimension, n=30}
\label{simulation_dimension_n=30}
\begin{center}
\begin{tabular}{|c|c|c|c|c|c|c|c|c|c|c|}
\hline
Case 1& 1 & 2 & 3 & 4 & {\bf 5} & 6 & 7 & 8 & 9 & 10 \\ \hline
$\lambda_1,\ldots,\lambda_{10}$ & 19 & 19 & 19 & 19 & {\bf 19} & 1 & 1 & 1 & 1 & 1 \\ \hline
C.1 \& $\hat{\bm{\tau}}^{(0)}$ & 0 & 0 & 104 & 9895 & {\bf 1} & 0 & 0 & 0 & 0 & 0 \\ \hline
C.1 \& $\hat{\bm{\tau}}^{(1)}$ & 0 & 0 & 0 & 9 & {\bf 9991} & 0 & 0 & 0 & 0 & 0 \\ \hline
C.1 \& $\hat{\bm{\tau}}^{(2)}$ & 0 & 0 & 0 & 339 & {\bf 9661} & 0 & 0 & 0 & 0 & 0 \\ \hline
C.2 \& $\hat{\bm{\tau}}^{(0)}$ & 0 & 0 & 285 & 6480 & {\bf 3235} & 0 & 0 & 0 & 0 & 0 \\ \hline
C.2 \& $\hat{\bm{\tau}}^{(1)}$ & 0 & 0 & 420 & 6346 & {\bf 3234} & 0 & 0 & 0 & 0 & 0 \\ \hline
C.2 \& $\hat{\bm{\tau}}^{(2)}$ & 0 & 0 & 2275 & 4930 & {\bf 2795} & 0 & 0 & 0 & 0 & 0 \\ \hline \hline
Case 2& 1 & 2 & 3 & {\bf 4} & 5 & 6 & 7 & 8 & 9 & 10 \\ \hline
$\lambda_1,\ldots,\lambda_{10}$  & 23.5 & 23.5 & 23.5 & {\bf 23.5} & 1 & 1 & 1 & 1 & 1 & 1 \\ \hline
C.1 \& $\hat{\bm{\tau}}^{(0)}$ & 0 & 0 & 7348 & {\bf 2652} & 0 & 0 & 0 & 0 & 0 & 0 \\ \hline
C.1 \& $\hat{\bm{\tau}}^{(1)}$ & 0 & 0 & 0 & {\bf 9828} & 172 & 0 & 0 & 0 & 0 & 0 \\ \hline
C.1 \& $\hat{\bm{\tau}}^{(2)}$ & 0 & 0 & 17 & {\bf 9982} & 1 & 0 & 0 & 0 & 0 & 0 \\ \hline
C.2 \& $\hat{\bm{\tau}}^{(0)}$ & 0 & 1 & 1087 & {\bf 8912} & 0 & 0 & 0 & 0 & 0 & 0 \\ \hline
C.2 \& $\hat{\bm{\tau}}^{(1)}$ & 0 & 5 & 1348 & {\bf 8647} & 0 & 0 & 0 & 0 & 0 & 0 \\ \hline
C.2 \& $\hat{\bm{\tau}}^{(2)}$ & 0 & 2 & 3801 & {\bf 6197} & 0 & 0 & 0 & 0 & 0 & 0 \\ \hline \hline
Case 3& 1 & 2 & {\bf 3} & 4 & 5 & 6 & 7 & 8 & 9 & 10 \\ \hline
$\lambda_1,\ldots,\lambda_{10}$  & 31 & 31 & {\bf 31} & 1 & 1 & 1 & 1 & 1 & 1 & 1 \\ \hline
C.1 \& $\hat{\bm{\tau}}^{(0)}$ & 0 & 376 & {\bf 9624} & 0 & 0 & 0 & 0 & 0 & 0 & 0 \\ \hline
C.1 \& $\hat{\bm{\tau}}^{(1)}$& 0 & 0 & {\bf 2410} & 6641 & 948 & 1 & 0 & 0 & 0 & 0 \\ \hline
C.1 \& $\hat{\bm{\tau}}^{(2)}$ & 0 & 0 & {\bf 9987} & 13 & 0 & 0 & 0 & 0 & 0 & 0 \\ \hline
C.2 \& $\hat{\bm{\tau}}^{(0)}$ & 0 & 22 & {\bf 9978} & 0 & 0 & 0 & 0 & 0 & 0 & 0 \\ \hline
C.2 \& $\hat{\bm{\tau}}^{(1)}$  & 0 & 50 & {\bf 9950} & 0 & 0 & 0 & 0 & 0 & 0 & 0 \\ \hline
C.2 \& $\hat{\bm{\tau}}^{(2)}$  & 0 & 30 & {\bf 9970} & 0 & 0 & 0 & 0 & 0 & 0 & 0 \\ \hline \hline
Case 4& 1 & {\bf 2} & 3 & 4 & 5 & 6 & 7 & 8 & 9 & 10 \\ \hline
$\lambda_1,\ldots,\lambda_{10}$  & 46 & {\bf 46} & 1 & 1 & 1 & 1 & 1 & 1 & 1 & 1 \\ \hline
C.1 \& $\hat{\bm{\tau}}^{(0)}$& 0 & {\bf 10000} & 0 & 0 & 0 & 0 & 0 & 0 & 0 & 0 \\ \hline
C.1 \& $\hat{\bm{\tau}}^{(1)}$ & 0 & {\bf 1} & 137 & 3545 & 5907 & 410 & 0 & 0 & 0 & 0 \\ \hline
C.1 \& $\hat{\bm{\tau}}^{(2)}$ & 0 & {\bf 8008} & 1903 & 88 & 1 & 0 & 0 & 0 & 0 & 0 \\ \hline
C.2 \& $\hat{\bm{\tau}}^{(0)}$ & 0 & {\bf 10000} & 0 & 0 & 0 & 0 & 0 & 0 & 0 & 0 \\ \hline
C.2 \& $\hat{\bm{\tau}}^{(1)}$  & 0 & {\bf 10000} & 0 & 0 & 0 & 0 & 0 & 0 & 0 & 0 \\ \hline
C.2 \& $\hat{\bm{\tau}}^{(2)}$  & 0 & {\bf 10000} & 0 & 0 & 0 & 0 & 0 & 0 & 0 & 0 \\ \hline \hline
Case 5& {\bf 1} & 2 & 3 & 4 & 5 & 6 & 7 & 8 & 9 & 10 \\ \hline
$\lambda_1,\ldots,\lambda_{10}$  & {\bf 91} & 1 & 1 & 1 & 1 & 1 & 1 & 1 & 1 & 1 \\ \hline
C.1 \& $\hat{\bm{\tau}}^{(0)}$ & {\bf 9991} & 9 & 0 & 0 & 0 & 0 & 0 & 0 & 0 & 0 \\ \hline
C.1 \& $\hat{\bm{\tau}}^{(1)}$ & {\bf 0} & 0 & 0 & 2 & 776 & 8170 & 1052 & 0 & 0 & 0 \\ \hline
C.1 \& $\hat{\bm{\tau}}^{(2)}$ & {\bf 1184} & 3743 & 3916 & 1064 & 93 & 0 & 0 & 0 & 0 & 0 \\ \hline
C.2 \& $\hat{\bm{\tau}}^{(0)}$ & {\bf 10000} & 0 & 0 & 0 & 0 & 0 & 0 & 0 & 0 & 0 \\ \hline
C.2 \& $\hat{\bm{\tau}}^{(1)}$  & {\bf 10000} & 0 & 0 & 0 & 0 & 0 & 0 & 0 & 0 & 0 \\ \hline
C.2 \& $\hat{\bm{\tau}}^{(2)}$  & {\bf 10000} & 0 & 0 & 0 & 0 & 0 & 0 & 0 & 0 & 0 \\ \hline \hline
\end{tabular}
\end{center}
\end{table}
\begin{table}
\begin{center}
\begin{tabular}{|c|c|c|c|c|c|c|c|c|c|c|}
\hline
Case 6& 1 & 2 & 3 & {\bf 4} & 5 & 6 & 7 & 8 & 9 & 10 \\ \hline
$\lambda_1,\ldots,\lambda_{10}$  & 27 & 23 & 19 & {\bf 15} & 11 & 1 & 1 & 1 & 1 & 1 \\ \hline
C.1 \& $\hat{\bm{\tau}}^{(0)}$ & 0 & 0 & 1018 & {\bf 8981} & 1 & 0 & 0 & 0 & 0 & 0 \\ \hline
C.1 \& $\hat{\bm{\tau}}^{(1)}$ & 0 & 0 & 0 & {\bf 48} & 9952 & 0 & 0 & 0 & 0 & 0 \\ \hline
C.1 \& $\hat{\bm{\tau}}^{(2)}$ & 0 & 0 & 0 & {\bf 1535} & 8465 & 0 & 0 & 0 & 0 & 0 \\ \hline
C.2 \& $\hat{\bm{\tau}}^{(0)}$ & 0 & 2 & 1483 & {\bf 7619} & 896 & 0 & 0 & 0 & 0 & 0 \\ \hline
C.2 \& $\hat{\bm{\tau}}^{(1)}$  & 0 & 12 & 1922 & {\bf 7171} & 895 & 0 & 0 & 0 & 0 & 0 \\ \hline
C.2 \& $\hat{\bm{\tau}}^{(2)}$  & 0 & 2 & 5407 &{\bf 3904} & 687 & 0 & 0 & 0 & 0 & 0 \\ \hline \hline
Case 7& 1 & 2 & 3 & {\bf 4} & 5 & 6 & 7 & 8 & 9 & 10 \\ \hline
$\lambda_1,\ldots,\lambda_{10}$  & 32 & 25.5 & 19 & {\bf 12.5} & 6 & 1 & 1 & 1 & 1 & 1 \\ \hline
C.1 \& $\hat{\bm{\tau}}^{(0)}$ & 0 & 0 & 5931 & {\bf 4069} & 0 & 0 & 0 & 0 & 0 & 0 \\ \hline
C.1 \& $\hat{\bm{\tau}}^{(1)}$ & 0 & 0 & 0 & {\bf 688} & 9311 & 1 & 0 & 0 & 0 & 0 \\ \hline
C.1 \& $\hat{\bm{\tau}}^{(2)}$ & 0 & 0 & 2 & {\bf 6764} & 3234 & 0 & 0 & 0 & 0 & 0 \\ \hline
C.2 \& $\hat{\bm{\tau}}^{(0)}$ & 0 & 35 & 4702 & {\bf 5254} & 9 & 0 & 0 & 0 & 0 & 0 \\ \hline
C.2 \& $\hat{\bm{\tau}}^{(1)}$  & 0 & 107 & 5209 & {\bf 4676} & 8 & 0 & 0 & 0 & 0 & 0 \\ \hline
C.2 \& $\hat{\bm{\tau}}^{(2)}$  & 0 & 49 & 8246 & {\bf 1701} & 4 & 0 & 0 & 0 & 0 & 0 \\ \hline \hline
Case 8& 1 & 2 & {\bf 3} & 4 & 5 & 6 & 7 & 8 & 9 & 10 \\ \hline
$\lambda_1,\ldots,\lambda_{10}$  & 37 & 28 & {\bf 19} & 10 & 1 & 1 & 1 & 1 & 1 & 1 \\ \hline
C.1 \& $\hat{\bm{\tau}}^{(0)}$ & 0 & 36 & {\bf 9847} & 117 & 0 & 0 & 0 & 0 & 0 & 0 \\ \hline
C.1 \& $\hat{\bm{\tau}}^{(1)}$ & 0 & 0 & {\bf 0} & 9013 & 987 & 0 & 0 & 0 & 0 & 0 \\ \hline
C.1 \& $\hat{\bm{\tau}}^{(2)}$ & 0 & 0 & {\bf 910} & 9090 & 0 & 0 & 0 & 0 & 0 & 0 \\ \hline
C.2 \& $\hat{\bm{\tau}}^{(0)}$ & 0 & 107 & {\bf 7496} & 2397 & 0 & 0 & 0 & 0 & 0 & 0 \\ \hline
C.2 \& $\hat{\bm{\tau}}^{(1)}$  & 0 & 317 & {\bf 7705} & 1978 & 0 & 0 & 0 & 0 & 0 & 0 \\ \hline
C.2 \& $\hat{\bm{\tau}}^{(2)}$  & 0 & 168 & {\bf 9347} & 485 & 0 & 0 & 0 & 0 & 0 & 0 \\ \hline \hline
Case 9& 1 & 2 & {\bf 3} & 4 & 5 & 6 & 7 & 8 & 9 & 10 \\ \hline
$\lambda_1,\ldots,\lambda_{10}$  & 46 & 31 & {\bf 16} & 1 & 1 & 1 & 1 & 1 & 1 & 1 \\ \hline
C.1 \& $\hat{\bm{\tau}}^{(0)}$ & 0 & 4471 & {\bf 5529} & 0 & 0 & 0 & 0 & 0 & 0 & 0 \\ \hline
C.1 \& $\hat{\bm{\tau}}^{(1)}$ & 0 & 0 & {\bf 598} & 6457 & 2923 & 22 & 0 & 0 & 0 & 0 \\ \hline
C.1 \& $\hat{\bm{\tau}}^{(2)}$ & 0 & 0 & {\bf 9971} & 29 & 0 & 0 & 0 & 0 & 0 & 0 \\ \hline
C.2 \& $\hat{\bm{\tau}}^{(0)}$ & 0 & 1022 & {\bf 8978} & 0 & 0 & 0 & 0 & 0 & 0 & 0 \\ \hline
C.2 \& $\hat{\bm{\tau}}^{(1)}$  & 0 & 1756 & {\bf 8244} & 0 & 0 & 0 & 0 & 0 & 0 & 0 \\ \hline
C.2 \& $\hat{\bm{\tau}}^{(2)}$  & 0 & 1234 & {\bf 8766} & 0 & 0 & 0 & 0 & 0 & 0 & 0 \\ \hline \hline
Case 10& {\bf 1} & 2 & 3 & 4 & 5 & 6 & 7 & 8 & 9 & 10 \\ \hline 
$\lambda_1,\ldots,\lambda_{10}$& {\bf 81} & 11 & 1 & 1 & 1 & 1 & 1 & 1 & 1 & 1 \\ \hline
C.1 \& $\hat{\bm{\tau}}^{(0)}$ & {\bf 6232} & 3767 & 1 & 0 & 0 & 0 & 0 & 0 & 0 & 0 \\ \hline
C.1 \& $\hat{\bm{\tau}}^{(1)}$ & {\bf 0} & 0 & 0 & 14 & 3144 & 6766 & 76 & 0 & 0 & 0 \\ \hline
C.1 \& $\hat{\bm{\tau}}^{(2)}$ & {\bf 0} & 4482 & 4664 & 824 & 30 & 0 & 0 & 0 & 0 & 0 \\ \hline
C.2 \& $\hat{\bm{\tau}}^{(0)}$ & {\bf 3939} & 6061 & 0 & 0 & 0 & 0 & 0 & 0 & 0 & 0 \\ \hline
C.2 \& $\hat{\bm{\tau}}^{(1)}$  & {\bf 5952} & 4048 & 0 & 0 & 0 & 0 & 0 & 0 & 0 & 0 \\ \hline
C.2 \& $\hat{\bm{\tau}}^{(2)}$  & {\bf 5171} & 4829 & 0 & 0 & 0 & 0 & 0 & 0 & 0 & 0 \\ \hline
\end{tabular}
\end{center}
\end{table}
\begin{table}
\caption{Histogram of Estimated Case, n=100}
\label{simulation_dimension_n=100}
\begin{center}
\begin{tabular}{|c|c|c|c|c|c|c|c|c|c|c|}
\hline
Case 1& 1 & 2 & 3 & 4 & {\bf 5} & 6 & 7 & 8 & 9 & 10 \\ \hline
$\lambda_1,\ldots,\lambda_{10}$ & 19 & 19 & 19 & 19 & {\bf 19} & 1 & 1 & 1 & 1 & 1 \\ \hline
C.1 \& $\hat{\bm{\tau}}^{(0)}$ & 0 & 0 & 0 & 9536 & {\bf 464} & 0 & 0 & 0 & 0 & 0 \\ \hline
C.1 \& $\hat{\bm{\tau}}^{(1)}$ & 0 & 0 & 0 & 1259 & {\bf 8741} & 0 & 0 & 0 & 0 & 0 \\ \hline
C.1 \& $\hat{\bm{\tau}}^{(2)}$ & 0 & 0 & 0 & 1922 & {\bf 8078} & 0 & 0 & 0 & 0 & 0 \\ \hline
C.2 \& $\hat{\bm{\tau}}^{(0)}$ & 0 & 0 & 0 & 50 & {\bf 9950} & 0 & 0 & 0 & 0 & 0 \\ \hline
C.2 \& $\hat{\bm{\tau}}^{(1)}$ & 0 & 0 & 0 & 50 & {\bf 9950} & 0 & 0 & 0 & 0 & 0 \\ \hline
C.2 \& $\hat{\bm{\tau}}^{(2)}$ & 0 & 0 & 0 & 50 & {\bf 9950} & 0 & 0 & 0 & 0 & 0 \\ \hline  \hline
Case 2& 1 & 2 & 3 & {\bf 4} & 5 & 6 & 7 & 8 & 9 & 10 \\ \hline 
$\lambda_1,\ldots,\lambda_{10}$ & 23.5 & 23.5 & 23.5 & {\bf 23.5} & 1 & 1 & 1 & 1 & 1 & 1 \\ \hline
C.1 \& $\hat{\bm{\tau}}^{(0)}$ & 0 & 0 & 287 & {\bf 9713} & 0 & 0 & 0 & 0 & 0 & 0 \\ \hline
C.1 \& $\hat{\bm{\tau}}^{(1)}$ & 0 & 0 & 0 & {\bf 10000} & 0 & 0 & 0 & 0 & 0 & 0 \\ \hline
C.1 \& $\hat{\bm{\tau}}^{(2)}$ & 0 & 0 & 9 & {\bf 9991} & 0 & 0 & 0 & 0 & 0 & 0 \\ \hline
C.2 \& $\hat{\bm{\tau}}^{(0)}$ & 0 & 0 & 1 & {\bf 9999} & 0 & 0 & 0 & 0 & 0 & 0 \\ \hline
C.2 \& $\hat{\bm{\tau}}^{(1)}$ & 0 & 0 & 1 & {\bf 9999} & 0 & 0 & 0 & 0 & 0 & 0 \\ \hline
C.2 \& $\hat{\bm{\tau}}^{(2)}$ & 0 & 0 & 2 & {\bf 9998} & 0 & 0 & 0 & 0 & 0 & 0 \\ \hline \hline
Case 3& 1 & 2 & {\bf 3} & 4 & 5 & 6 & 7 & 8 & 9 & 10 \\ \hline
$\lambda_1,\ldots,\lambda_{10}$ & 31 & 31 & {\bf 31} & 1 & 1 & 1 & 1 & 1 & 1 & 1 \\ \hline
C.1 \& $\hat{\bm{\tau}}^{(0)}$ & 0 & 0 & {\bf 10000} & 0 & 0 & 0 & 0 & 0 & 0 & 0 \\ \hline
C.1 \& $\hat{\bm{\tau}}^{(1)}$ & 0 & 0 & {\bf 10000} & 0 & 0 & 0 & 0 & 0 & 0 & 0 \\ \hline
C.1 \& $\hat{\bm{\tau}}^{(2)}$ & 0 & 0 & {\bf 10000} & 0 & 0 & 0 & 0 & 0 & 0 & 0 \\ \hline
C.2 \& $\hat{\bm{\tau}}^{(0)}$ & 0 & 0 & {\bf 10000} & 0 & 0 & 0 & 0 & 0 & 0 & 0 \\ \hline
C.2 \& $\hat{\bm{\tau}}^{(1)}$ & 0 & 0 & {\bf 10000} & 0 & 0 & 0 & 0 & 0 & 0 & 0 \\ \hline
C.2 \& $\hat{\bm{\tau}}^{(2)}$ & 0 & 0 & {\bf 10000} & 0 & 0 & 0 & 0 & 0 & 0 & 0 \\ \hline \hline
Case 4& 1 & {\bf 2} & 3 & 4 & 5 & 6 & 7 & 8 & 9 & 10 \\ \hline
$\lambda_1,\ldots,\lambda_{10}$ & 46 & {\bf 46} & 1 & 1 & 1 & 1 & 1 & 1 & 1 & 1 \\ \hline
C.1 \& $\hat{\bm{\tau}}^{(0)}$ & 0 & {\bf 10000} & 0 & 0 & 0 & 0 & 0 & 0 & 0 & 0 \\ \hline
C.1 \& $\hat{\bm{\tau}}^{(1)}$ & 0 & {\bf 10000} & 0 & 0 & 0 & 0 & 0 & 0 & 0 & 0 \\ \hline
C.1 \& $\hat{\bm{\tau}}^{(2)}$ & 0 & {\bf 10000} & 0 & 0 & 0 & 0 & 0 & 0 & 0 & 0 \\ \hline
C.2 \& $\hat{\bm{\tau}}^{(0)}$ & 0 & {\bf 10000} & 0 & 0 & 0 & 0 & 0 & 0 & 0 & 0 \\ \hline
C.2 \& $\hat{\bm{\tau}}^{(1)}$ & 0 & {\bf 10000} & 0 & 0 & 0 & 0 & 0 & 0 & 0 & 0 \\ \hline
C.2 \& $\hat{\bm{\tau}}^{(2)}$ & 0 & {\bf 10000} & 0 & 0 & 0 & 0 & 0 & 0 & 0 & 0 \\ \hline \hline
Case 5 & {\bf 1} & 2 & 3 & 4 & 5 & 6 & 7 & 8 & 9 & 10 \\ \hline
$\lambda_1,\ldots,\lambda_{10}$ & {\bf 91} & 1 & 1 & 1 & 1 & 1 & 1 & 1 & 1 & 1 \\ \hline
C.1 \& $\hat{\bm{\tau}}^{(0)}$ & {\bf 10000} & 0 & 0 & 0 & 0 & 0 & 0 & 0 & 0 & 0 \\ \hline
C.1 \& $\hat{\bm{\tau}}^{(1)}$ & {\bf 9995} & 5 & 0 & 0 & 0 & 0 & 0 & 0 & 0 & 0 \\ \hline
C.1 \& $\hat{\bm{\tau}}^{(2)}$ & {\bf 9999} & 1 & 0 & 0 & 0 & 0 & 0 & 0 & 0 & 0 \\ \hline
C.2 \& $\hat{\bm{\tau}}^{(0)}$ & {\bf 10000} & 0 & 0 & 0 & 0 & 0 & 0 & 0 & 0 & 0 \\ \hline
C.2 \& $\hat{\bm{\tau}}^{(1)}$ & {\bf 10000} & 0 & 0 & 0 & 0 & 0 & 0 & 0 & 0 & 0 \\ \hline
C.2 \& $\hat{\bm{\tau}}^{(2)}$ & {\bf 10000} & 0 & 0 & 0 & 0 & 0 & 0 & 0 & 0 & 0 \\ \hline
\end{tabular}
\end{center}
\end{table}
\begin{table}
\begin{center}
\begin{tabular}{|c|c|c|c|c|c|c|c|c|c|c|}
\hline
Case 6 & 1 & 2 & 3 & {\bf 4} & 5 & 6 & 7 & 8 & 9 & 10 \\ \hline
$\lambda_1,\ldots,\lambda_{10}$ & 27 & 23 & 19 & {\bf 15} & 11 & 1 & 1 & 1 & 1 & 1 \\ \hline
C.1 \& $\hat{\bm{\tau}}^{(0)}$ & 0 & 0 & 0 & {\bf 9998} & 2 & 0 & 0 & 0 & 0 & 0 \\ \hline
C.1 \& $\hat{\bm{\tau}}^{(1)}$ & 0 & 0 & 0 & {\bf 8849} & 1151 & 0 & 0 & 0 & 0 & 0 \\ \hline
C.1 \& $\hat{\bm{\tau}}^{(2)}$ & 0 & 0 & 0 & {\bf 9422} & 578 & 0 & 0 & 0 & 0 & 0 \\ \hline
C.2 \& $\hat{\bm{\tau}}^{(0)}$ & 0 & 0 & 15 & {\bf 5304} & 4681 & 0 & 0 & 0 & 0 & 0 \\ \hline
C.2 \& $\hat{\bm{\tau}}^{(1)}$ & 0 & 0 & 20 & {\bf 5299} & 4681 & 0 & 0 & 0 & 0 & 0 \\ \hline
C.2 \& $\hat{\bm{\tau}}^{(2)}$ & 0 & 0 & 87 & {\bf 5239} & 4674 & 0 & 0 & 0 & 0 & 0 \\ \hline \hline
Case 7& 1 & 2 & 3 & {\bf 4} & 5 & 6 & 7 & 8 & 9 & 10 \\ \hline
$\lambda_1,\ldots,\lambda_{10}$ & 32 & 25.5 & 19 & {\bf 12.5} & 6 & 1 & 1 & 1 & 1 & 1 \\ \hline
C.1 \& $\hat{\bm{\tau}}^{(0)}$ & 0 & 0 & 1113 & {\bf 8887} & 0 & 0 & 0 & 0 & 0 & 0 \\ \hline
C.1 \& $\hat{\bm{\tau}}^{(1)}$ & 0 & 0 & 1 & {\bf 9999} & 0 & 0 & 0 & 0 & 0 & 0 \\ \hline
C.1 \& $\hat{\bm{\tau}}^{(2)}$ & 0 & 0 & 25 & {\bf 9975} & 0 & 0 & 0 & 0 & 0 & 0 \\ \hline
C.2 \& $\hat{\bm{\tau}}^{(0)}$ & 0 & 0 & 1095 & {\bf 8905} & 0 & 0 & 0 & 0 & 0 & 0 \\ \hline
C.2 \& $\hat{\bm{\tau}}^{(1)}$ & 0 & 0 & 1228 & {\bf 8772} & 0 & 0 & 0 & 0 & 0 & 0 \\ \hline
C.2 \& $\hat{\bm{\tau}}^{(2)}$ & 0 & 0 & 2321 & {\bf 7679} & 0 & 0 & 0 & 0 & 0 & 0 \\ \hline \hline
Case 8& 1 &  2 & {\bf 3} & 4 & 5 & 6 & 7 & 8 & 9 & 10 \\ \hline
$\lambda_1,\ldots,\lambda_{10}$ & 37 & 28 & {\bf 19} & 10 & 1 & 1 & 1 & 1 & 1 & 1 \\ \hline
C.1 \& $\hat{\bm{\tau}}^{(0)}$ & 0 & 0 & {\bf 9964} & 36 & 0 & 0 & 0 & 0 & 0 & 0 \\ \hline
C.1 \& $\hat{\bm{\tau}}^{(1)}$ & 0 & 0 & {\bf 5754} & 4246 & 0 & 0 & 0 & 0 & 0 & 0 \\ \hline
C.1 \& $\hat{\bm{\tau}}^{(2)}$ & 0 & 0 & {\bf 8784} & 1216 & 0 & 0 & 0 & 0 & 0 & 0 \\ \hline
C.2 \& $\hat{\bm{\tau}}^{(0)}$ & 0 & 0 & {\bf 6280} & 3720 & 0 & 0 & 0 & 0 & 0 & 0 \\ \hline
C.2 \& $\hat{\bm{\tau}}^{(1)}$ & 0 & 0 & {\bf 6531} & 3469 & 0 & 0 & 0 & 0 & 0 & 0 \\ \hline
C.2 \& $\hat{\bm{\tau}}^{(2)}$ & 0 & 0 & {\bf 7859} & 2141 & 0 & 0 & 0 & 0 & 0 & 0 \\ \hline \hline
Case 9& 1 & 2 & {\bf 3} & 4 & 5 & 6 & 7 & 8 & 9 & 10 \\ \hline
$\lambda_1,\ldots,\lambda_{10}$ & 46 & 31 & {\bf 16} & 1 & 1 & 1 & 1 & 1 & 1 & 1 \\ \hline
C.1 \& $\hat{\bm{\tau}}^{(0)}$ & 0 & 1628 &{\bf 8372} & 0 & 0 & 0 & 0 & 0 & 0 & 0 \\ \hline
C.1 \& $\hat{\bm{\tau}}^{(1)}$ & 0 & 0 & {\bf 10000} & 0 & 0 & 0 & 0 & 0 & 0 & 0 \\ \hline
C.1 \& $\hat{\bm{\tau}}^{(2)}$ & 0 & 49 & {\bf 9951} & 0 & 0 & 0 & 0 & 0 & 0 & 0 \\ \hline
C.2 \& $\hat{\bm{\tau}}^{(0)}$ & 0 & 19 & {\bf 9981} & 0 & 0 & 0 & 0 & 0 & 0 & 0 \\ \hline
C.2 \& $\hat{\bm{\tau}}^{(1)}$ & 0 & 38 & {\bf 9962} & 0 & 0 & 0 & 0 & 0 & 0 & 0 \\ \hline
C.2 \& $\hat{\bm{\tau}}^{(2)}$ & 0 & 22 & {\bf 9978} & 0 & 0 & 0 & 0 & 0 & 0 & 0 \\ \hline \hline
Case 10& {\bf 1} & 2 & 3 & 4 & 5 & 6 & 7 & 8 & 9 & 10 \\ \hline
$\lambda_1,\ldots,\lambda_{10}$ & {\bf 81} & 11 & 1 & 1 & 1 & 1 & 1 & 1 & 1 & 1 \\ \hline
C.1 \& $\hat{\bm{\tau}}^{(0)}$ & {\bf 6679} & 3321 & 0 & 0 & 0 & 0 & 0 & 0 & 0 & 0 \\ \hline
C.1 \& $\hat{\bm{\tau}}^{(1)}$ & {\bf 208} & 9792 & 0 & 0 & 0 & 0 & 0 & 0 & 0 & 0 \\ \hline
C.1 \& $\hat{\bm{\tau}}^{(2)}$ & {\bf 1108} & 8892 & 0 & 0 & 0 & 0 & 0 & 0 & 0 & 0 \\ \hline
C.2 \& $\hat{\bm{\tau}}^{(0)}$ & {\bf 2952} & 7048 & 0 & 0 & 0 & 0 & 0 & 0 & 0 & 0 \\ \hline
C.2 \& $\hat{\bm{\tau}}^{(1)}$ & {\bf 4039} & 5961 & 0 & 0 & 0 & 0 & 0 & 0 & 0 & 0 \\ \hline
C.2 \& $\hat{\bm{\tau}}^{(2)}$ & {\bf 3612} & 6388 & 0 & 0 & 0 & 0 & 0 & 0 & 0 & 0 \\ \hline 
\end{tabular}
\end{center}
\end{table}
Table \ref{simulation_dimension_n=30} and \ref{simulation_dimension_n=100} are the result of the simulation. To explain the meaning of  each number in the table, take the first case in  Table \ref{simulation_dimension_n=30} as an example, where $\lambda_i,\ i=1,\ldots,10$ are given by $19, 19, 19, 19, 19, 1, 1, 1, 1, 1$. (Note that for every case, $\lambda_i$'s are designed so that $\sum_{i=1}^p \lambda_i =100$, hence each $\lambda_i$ equals $\tau_i$ in percent figures.) The row of ``C.1 \& $\hat{\bm{\tau}}^{(0)}$'' is the histogram of the estimated dimension by the combination of criterion 1 and $\hat{\bm{\tau}}^{(0)}$. The boldface position ( in this case ``{\bf 1}'' ) indicates that the ``true'' dimension decided by criterion 1 from the population eigenvalues, $\lambda_i$'s. In this simulation,  the ``true'' dimensions are designed to take the same value by either criterion 1 or 2.

We can observe following points from Table  \ref{simulation_dimension_n=30}. With respect to the criterion 1, we notice that the classical estimator, $\hat{\bm{\tau}}^{(0)}$, tends to underestimate the dimension (see Case 1, 2, 7, 9), while $\hat{\bm{\tau}}^{(1)}$ tend to overestimate it  (see Case 3--10). From \eqref{relation_tau1_tau2}, we notice that $\hat{\bm{\tau}}^{(2)}$ is located between $\hat{\bm{\tau}}^{(0)}$ and $\hat{\bm{\tau}}^{(1)}.$ Though  $\hat{\bm{\tau}}^{(2)}$ is still likely to overestimate the dimension (see Case 5, 6, 8, 10), the tendency is weakened compared to  $\hat{\bm{\tau}}^{(1)}$. 
On the criterion 2, we can not find as significant a difference as criterion 1 among three estimators. Every estimator tends to underestimate the dimension in some cases (see Case 1, 7) and overestimate it in other cases (see Case 10).  

In most cases in Table \ref{simulation_dimension_n=100}, the estimation for dimension is made correctly. However, despite a high degree of freedom, we still observe the tendency of  $\hat{\bm{\tau}}^{(0)}$ to underestimation (see Case 1) and that of  $\hat{\bm{\tau}}^{(1)}$ or  $\hat{\bm{\tau}}^{(2)}$ to overestimation (see Case 10) with respect to the criterion 1.

 Though both underestimation and overestimation are undesirable, the former is more crucial, and have more substantial effect on the results obtained in  principal component analysis or factor analysis, since important component (factor) is neglected  (see e.g. the comment in p 278 in Fabriger et. al. (1999)). In 
this sense,  $\hat{\bm{\tau}}^{(1)}$ and  $\hat{\bm{\tau}}^{(2)}$ are superior to  $\hat{\bm{\tau}}^{(0)}$. The tendency to overestimation of  $\hat{\bm{\tau}}^{(1)}$ is weakened in  $\hat{\bm{\tau}}^{(2)}$ since $\beta_i^{(2)}$'s are closer to unit.  We can correct the overestimation further by selecting $\beta_i$'s that are much closer to unit, but still satisfy the three conditions in Theorem \ref{dominance_result}. 
\section{Conclusion}
We can summarize the results of this paper as follows;
\begin{enumerate}
\item The distribution of the sample contribution rates is identical within a family of elliptically contoured distributions. It is determined solely by the population contribution rates.
\item A class of new estimators of the population contribution rates was derived. In the estimation of the normalized population covariance matrix, the estimator composed of the new estimator and the sample eigenvectors dominates the estimator composed of the classical estimator and the sample eigenvectors with respect to the entropy loss function.
\item A simulation study shows that the new estimators perform substantially better than the classical estimator with respect to the risk derived from a quadratic loss function. Another simulation study shows that the new estimators tend to overestimate the dimension. They are more suitable than the classical estimator for the decision of dimension in principal component analysis or factor analysis, since the classical estimator is likely to underestimate the dimension.
\end{enumerate}
\section{Acknowledgment}
We really appreciate many constructive comments and suggestions from two anonymous referees. They broadened the author's vision and improved the quality of this paper.
\section{Appendix}
\subsection{Proof of \eqref{assym_expan_E(d_i)}}
We can suppose $\sss\sim W_p(n,\bm{\Lambda}),\ \bm{\Lambda}={\rm diag}(\lambda_1,\ldots,\lambda_p)$, since the distribution of $d_i$'s depends only on $\bm{\lambda}$. Let $\bm{\Delta}=(\delta_{ij})$ be defined as
$$
\bm{\Delta}=n^{-1}\sss-\bm{\Lambda}=\bm{A}-\bm{\Lambda}.
$$
(3) of Lawley (1959) gives the following expansion of $l_r^*\ (1\leq r \leq p)$,
\begin{equation}
\label{expansion_l}
\begin{split}
l_r^*&=\lambda_r+\delta_{rr}+\sum_{i\ne r}\frac{\delta_{ri}^2}{\lambda_r-\lambda_i}\\
&\quad-\delta_{rr}\sum_{i\ne r}\Bigl(\frac{\delta_{ri}}{\lambda_r-\lambda_i}\Bigr)^2+\sum_{i\ne r}\sum_{j \ne r}\Bigl(\frac{\delta_{ri}\delta_{rj}\delta_{ij}}{(\lambda_r-\lambda_i)(\lambda_r-\lambda_j)}\Bigr)+O(\norm{\bm{\delta}}^4).
\end{split}
\end{equation}
Since $\sum_{i=1}^p l_i^*=\sum_{i=1}^p(\lambda_i+\delta_{ii})$, we have the Taylor expansion of $(\sum_{i=1}^pl_i^*)^{-1}$ 
\begin{equation}
\label{expansion_trl_inverse}
\begin{split}
\frac{1}{\sum_{i=1}^pl_i^*} &=\frac{1}{\sum_{i=1}^p(\lambda_i+\delta_{ii})}\\
&=\frac{1}{\sum_{i=1}^p \lambda_i}
-\frac{\sum_{i=1}^p\delta_{ii}}{(\sum_{i=1}^p \lambda_i)^2}+2\frac{\sum_{_i<j}\delta_{ii}\delta_{jj}}{(\sum_{i=1}^p \lambda_i)^3}\\
&\qquad+\frac{\sum_{i=1}^p\delta_{ii}^2}{(\sum_{i=1}^p \lambda_i)^3}-\frac{\sum_{1\leq i,j,k \leq p}\delta_{ii}\delta_{jj}\delta_{kk}}{(\sum_{i=1}^p \lambda_i)^4}+O(\norm{\bm{\delta}}^4).
\end{split}
\end{equation}
Combining \eqref{expansion_l} and \eqref{expansion_trl_inverse}, we have 
\begin{equation}
\label{expansion_contri_rate}
\begin{split}
d_r&=\frac{l_r^*}{\sum_{i=1}^p l_i^*}\\&=\frac{\lambda_r}{\sum_{i=1}^p\lambda_i}+\frac{1}{\sum_{i=1}^p\lambda_i}\delta_{rr}-\frac{\lambda_r}{(\sum_{i=1}^p\lambda_i)^2}\sum_{i=1}^p\delta_{ii}\\
&\quad +\frac{\lambda_r}{(\sum_{i=1}^p\lambda_i)^3}\Bigl(\sum_{i=1}^p \delta_{ii}^2+2\sum_{i<j}\delta_{ii}\delta_{jj}\Bigr)\\
&\quad -\frac{1}{(\sum_{i=1}^p\lambda_i)^2}\Bigl(\sum_{i=1}^p\delta_{ii}\Bigr)\delta_{rr}+\frac{1}{\sum_{i=1}^p\lambda_i}\sum_{i\ne r}\frac{\delta_{ri}^2}{\lambda_r-\lambda_i}+O(\norm{\bm{\delta}}^3).
\end{split}
\end{equation}
We can easily calculate the low-dimensional moments of $\bm{\delta}$. They are given as follows;
\begin{equation}
\label{moment_delta_1}
E[\delta_{ii}]=0,\quad 1\leq i \leq p.
\end{equation}
\begin{equation}
\label{moment_delta_2}
E[\delta_{ii}\delta_{jj}]=\left\{
\begin{split}
&2n^{-1}\lambda_{i}^2,\quad &\text{if $i=j$,}\\
&0,\quad &\text{if $i\ne j$.}
\end{split}
\right.
\end{equation}
\begin{equation}
\label{moment_delta_3}
E[\delta_{ij}^2]=n^{-1} \lambda_i \lambda_j, \quad \text{if $i\ne j.$}
\end{equation}
In addition, 
\begin{equation}
\label{moment_delta_k}
E[\norm{\bm{\delta}}^k]=O(n^{-(k-1)}).
\end{equation}
Substituting \eqref{moment_delta_1} -- \eqref{moment_delta_k} into \eqref{expansion_contri_rate}, we have the desired result.  
\hfill\rule{5pt}{10pt}
\subsection{Proof of \eqref{density_S}}
Let the singular value decomposition of $\bm{Z}$ denoted by
$$
\bm{Z}=\bm{O}\bm{D}\bm{H},\quad \bm{D}={\rm diag}(d_1,\ldots,d_p),\ d_1>\cdots>d_p>0,\quad \bm{O}\in V_{p,n},\ \bm{H}\in O(p),
$$
where $O(p)$ is the set of  $p$-dimensional orthogonal matrices. From the Jacobian of this decomposition (see e.g. Theorem 5 of Uhlig (1994)), we have
$$
d\bm{Z} \propto \prod_{i=1}^p d_i^{n-p} \prod_{i<j} (d_i^2 -d _j^2)\: \mu_{p,p}(d\bm{H})\: \mu_{p,n}(d\bm{O})\: d\bm{d},
$$
where 
$\mu_{p,p}$ and $\mu_{p,n}$ are the invariant probability measures respectively on $O(p)$ and $V_{p,n}$. 
Further by the transformation $d_i \to t_i=d_i^2\ i=1,\ldots,p$, we have
$$
d\bm{Z} \propto \prod_{i=1}^p t_i^{(n-p-1)/2} \prod_{i<j}(t_i-t_j) \: \mu_{p,p}(d\bm{H})\: \mu_{p,n}(d\bm{O})\: d\bm{t}.
$$
Notice that
\begin{equation}
\label{spectral_decomp}
\sss=\bm{H}'\bm{T}\bm{H},\quad \bm{T}={\rm diag}(t_1,\ldots,t_p), 
\end{equation}
hence, 
\begin{align*}
&f({\rm tr}\bm{Z}'\bm{Z}\s^{-1})|\s|^{-n/2}\:d\bm{Z}\\
&\propto  f({\rm tr}\sss\s^{-1})|\sss|^{(n-p-1)/2}|\s|^{-n/2}\prod_{i<j}(t_i-t_j) \: \mu_{p,p}(d\bm{H})\: \mu_{p,n}(d\bm{O})\: d\bm{t}.
\end{align*}
If we integrate the right side of the above equation over $V_{p,n}$, we have the following density function of $\bm{H}$ and $\bm{t}$ with respect to $\mu_{p,p}(d\bm{H})d\bm{t}$
\begin{equation}
\label{density_H_t}
c_0 f({\rm tr}\sss\s^{-1})|\sss|^{(n-p-1)/2}|\s|^{-n/2}\prod_{i<j} (t_i-t_j)
\end{equation}
with some constant $c_0.$
Combined with a formula about the spectral decomposition \eqref{spectral_decomp} (see e.g. (22) on the p105 of Muirhead(1982))
$$
d\sss \propto \prod_{i<j}(t_i-t_j)\: \mu_{pp}(d\bm{H})\:d\bm{t},
$$
\eqref{density_H_t} leads to the following density of $\sss$ with respect to Lebesgue measure
\begin{equation}
c_1 f({\rm tr}\sss\s^{-1}) |\sss|^{(n-p-1)/2}|\s|^{-n/2}
\end{equation}
with a normalizing constant $c_1.$ 
\hfill \rule{5pt}{10pt}

\end{document}